\documentclass[11pt,fleqn]{article}
\usepackage{amsfonts, epsfig}
\newcommand{\ol}{\setlength{\itemsep}{0pt.}\begin{enumerate}}
\newcommand{\eol}{\end{enumerate}\setlength{\itemsep}{-\parsep}}
\newcommand{\ignore}[1]{}
\title{An inequality for functions on the Hamming cube}
\author{Alex Samorodnitsky
\thanks{School of Engineering and Computer Science,
The Hebrew University of Jerusalem,
Jerusalem 91904, Israel.
Research partially supported by ISF
grant 1241/11 and by BSF grant 2010451. Part of this work was done while the author was a fellow at the Radcliffe Institute in 2009-2010. }
}

\begin{document}
\date{}
\maketitle


\newtheorem{THEOREM}{Theorem}[section]
\newenvironment{theorem}{\begin{THEOREM} \hspace{-.85em} {\bf :}
}%
                        {\end{THEOREM}}
\newtheorem{LEMMA}[THEOREM]{Lemma}
\newenvironment{lemma}{\begin{LEMMA} \hspace{-.85em} {\bf :} }%
                      {\end{LEMMA}}
\newtheorem{COROLLARY}[THEOREM]{Corollary}
\newenvironment{corollary}{\begin{COROLLARY} \hspace{-.85em} {\bf
:} }%
                          {\end{COROLLARY}}
\newtheorem{PROPOSITION}[THEOREM]{Proposition}
\newenvironment{proposition}{\begin{PROPOSITION} \hspace{-.85em}
{\bf :} }%
                            {\end{PROPOSITION}}
\newtheorem{DEFINITION}[THEOREM]{Definition}
\newenvironment{definition}{\begin{DEFINITION} \hspace{-.85em} {\bf
:} \rm}%
                            {\end{DEFINITION}}
\newtheorem{EXAMPLE}[THEOREM]{Example}
\newenvironment{example}{\begin{EXAMPLE} \hspace{-.85em} {\bf :}
\rm}%
                            {\end{EXAMPLE}}
\newtheorem{CONJECTURE}[THEOREM]{Conjecture}
\newenvironment{conjecture}{\begin{CONJECTURE} \hspace{-.85em}
{\bf :} \rm}%
                            {\end{CONJECTURE}}
\newtheorem{MAINCONJECTURE}[THEOREM]{Main Conjecture}
\newenvironment{mainconjecture}{\begin{MAINCONJECTURE} \hspace{-.85em}
{\bf :} \rm}%
                            {\end{MAINCONJECTURE}}
\newtheorem{PROBLEM}[THEOREM]{Problem}
\newenvironment{problem}{\begin{PROBLEM} \hspace{-.85em} {\bf :}
\rm}%
                            {\end{PROBLEM}}
\newtheorem{QUESTION}[THEOREM]{Question}
\newenvironment{question}{\begin{QUESTION} \hspace{-.85em} {\bf :}
\rm}%
                            {\end{QUESTION}}
\newtheorem{REMARK}[THEOREM]{Remark}
\newenvironment{remark}{\begin{REMARK} \hspace{-.85em} {\bf :}
\rm}%
                            {\end{REMARK}}

\newcommand{\thm}{\begin{theorem}}
\newcommand{\lem}{\begin{lemma}}
\newcommand{\pro}{\begin{proposition}}
\newcommand{\dfn}{\begin{definition}}
\newcommand{\rem}{\begin{remark}}
\newcommand{\xam}{\begin{example}}
\newcommand{\cnj}{\begin{conjecture}}
\newcommand{\mcnj}{\begin{mainconjecture}}
\newcommand{\prb}{\begin{problem}}
\newcommand{\que}{\begin{question}}
\newcommand{\cor}{\begin{corollary}}
\newcommand{\prf}{\noindent{\bf Proof:} }
\newcommand{\ethm}{\end{theorem}}
\newcommand{\elem}{\end{lemma}}
\newcommand{\epro}{\end{proposition}}
\newcommand{\edfn}{\bbox\end{definition}}
\newcommand{\erem}{\bbox\end{remark}}
\newcommand{\exam}{\bbox\end{example}}
\newcommand{\ecnj}{\bbox\end{conjecture}}
\newcommand{\emcnj}{\bbox\end{mainconjecture}}
\newcommand{\eprb}{\bbox\end{problem}}
\newcommand{\eque}{\bbox\end{question}}
\newcommand{\ecor}{\end{corollary}}
\newcommand{\eprf}{\bbox}
\newcommand{\beqn}{\begin{equation}}
\newcommand{\eeqn}{\end{equation}}
\newcommand{\wbox}{\mbox{$\sqcap$\llap{$\sqcup$}}}
\newcommand{\bbox}{\vrule height7pt width4pt depth1pt}
\newcommand{\qed}{\bbox}
\def\sup{^}

\def\H{\{0,1\}^n}
\def\B{\{0,1\}}

\def\S{S(n,w)}

\def\n{\lfloor \frac n2 \rfloor}

\def \E{{\cal E}}
\def \Ex{\mathbb E}
\def \R{\mathbb R}
\def \Z{\mathbb Z}
\def \F{\mathbb F}

\def \noi{\noindent}

\def\<{\left<}
\def\>{\right>}
\def \({\left(}
\def \){\right)}
\def \e{\epsilon}
\def \r{\rfloor}

\def \1{{\bf 1}}

\def\Tp{Tchebyshef polynomial}
\def\Tps{TchebysDeto be the maximafine $A(n,d)$ l size of a code with distance $d$hef polynomials}
\newcommand{\rarrow}{\rightarrow}

\newcommand{\larrow}{\leftarrow}

\overfullrule=0pt
\def\setof#1{\lbrace #1 \rbrace}

\begin{abstract}
We prove an inequality for functions on the discrete cube $\H$ extending the edge-isoperimetric inequality for sets.

This inequality turns out to be equivalent to the following claim about random walks on the cube: Subcubes maximize 'mean first exit time' among all subsets of the cube of the same cardinality.
\end{abstract}

\section{Introduction}

\noi Isoperimetric inequalities play an important role in describing the geometry of ambient spaces \cite{Bezrukov, Ros}. This paper deals with one such space, the discrete cube $\H$.  This is a graph with $2^n$ vertices indexed by boolean strings of length $n$. Two vertices are connected by an edge if they differ in one coordinate. The {\it edge-isoperimetric} inequality \cite{Harper} for $\H$ provides well-known example for a discrete isoperimetric inequality.

\noi The edge boundary $\partial A$ of a subset $A \subseteq \H$ is the set of edges between $A$ and its complement. The edge-isoperimetric inequality compares between the cardinality of the set and of its boundary:
\beqn
\label{ineq:edge-isop}
|\partial A| \ge |A| \cdot \log_2\(\frac{2^n}{|A|}\)
\eeqn
One of its implications is that a simple random walk in the cube doesn't stay for too long in any given subset. This can be used to prove upper bounds on the mixing time of the walk \cite{Jerrum}.

\noi This inequality can also be viewed as an inequality for characteristic functions on $\H$. For a function $g:~\H \rarrow \R$, let the  {\it Dirichlet quadratic form} of $g$ be given by
$$
{\cal E}(g,g) = \Ex_x \sum_{y\sim x} (g(x) - g(y))^2
$$
Here the expectation is taken with respect to the uniform probability measure on the cube. The notation $x \sim y$ means that $x$ and $y$ are connected by an edge. Then (\ref{ineq:edge-isop}) can be rewritten, for $g = 1_A$ as
\beqn
\label{ineq:weak-LS-2}
{\cal E}(g,g) \ge 2 \Ex g^2 \cdot \log_2\(\frac{\Ex g^2}{\Ex^2 |g|}\)
\eeqn

\noi It is natural to look for inequalities for real-valued functions $g$ on the cube generalizing (\ref{ineq:edge-isop}). One such inequality is the logarithmic Sobolev inequality \cite{Gross} :
\[
{\cal E}(g,g) \ge Ent\(g^2\) = \Ex g^2 \ln g^2 - \Ex g^2 \ln \Ex g^2
\]
For $g = 1_A$ this becomes $|\partial A| \ge |A| \cdot \ln\(\frac{2^n}{|A|}\)$, recovering (\ref{ineq:edge-isop}) up to a multiplicative factor of $1/\ln 2$.

\noi For a general real-valued function $g$, the logarithmic Sobolev inequality has been observed \cite{FS, Montenegro-Tetali} to imply
\beqn
\label{ineq:weak-LS}
{\cal E}(g,g) \ge 2  \Ex g^2 \cdot \ln\(\frac{\Ex g^2}{\Ex^2 |g|}\)
\eeqn
This extends (\ref{ineq:weak-LS-2}), again up to a multiplicative factor of $1/\ln 2$.

\noi It is useful to look for inequalities for general functions reducing to an isoperimetric inequality with the {\it correct} constant in the special case of characteristic functions. Such an inequality would, in particular, mean that the characteristic function of an {\it isoperimetric} set\footnote{Recall that a set is isoperimetric if it satisfies an isoperimetric inequality with equality.}, or an "almost-isoperimetric" set, is an optimal (or nearly optimal) solution of a continuous extremal problem, and as such, might be expected to have an interesting structure. We refer to \cite{Barthe} for an example of relevant work in continuous analysis.

\noi As observed in \cite{FS}, the inequality (\ref{ineq:weak-LS}) is in fact tight for general real-valued functions. Therefore, to recover correct constants, we need to look for different extensions of (\ref{ineq:edge-isop}).

\noi This paper gives one example of such a an inequality.

\thm
\label{thm:isop-ineq}
Let $A$ be a subset of $\H$ and let $g$ be a real-valued function on $\H$ supported on $A$. Then
\beqn
\label{ineq:zero-outside}
{\cal E}(g,g) \ge 2 \cdot \frac{1}{2^n \cdot |A|} \log_2\(\frac{2^{n}}{|A|}\) \cdot \(\sum_{x \in A} |g(x)|\)^2
\eeqn
\ethm

\noi The dependence on $g$ on the right hand side of this inequality is weaker than that in the logarithmic Sobolev inequality, or that in (\ref{ineq:weak-LS}). However, it does give the right constant. In fact, substituting $g = 1_A$ recovers (\ref{ineq:edge-isop}).

\noi It turns out that (\ref{ineq:zero-outside}) is equivalent to a statement about random walks in the cube. Let $A$ be a subset of $\H$. Let $Y$ be a random variable defined as follows: choose a uniformly random point $a \in A$ and consider the simple random walk in $\H$ starting from $a$. Then $Y$ measures the time it takes the walk to exit $A$ for the first time. We refer to $\Ex Y$ as the {\it mean first exit time} of $A$. This is a parameter of a subset $A$ of the cube.

The following claim is equivalent to Theorem~\ref{thm:isop-ineq}.
\thm
\label{thm:random-walk}
Subcubes maximize mean first exit time among all subsets of the cube of the same cardinality.

More precisely, for any subset $A$ of $\H$,
\beqn
\label{ineq:random-walk}
\Ex Y \le \frac{n}{\log_2\(\frac{2^n}{|A|}\)}
\eeqn
If $A$ is a subcube, this is an equality.
\ethm

\noi This paper is organized as follows. We show equivalence of theorems~\ref{thm:isop-ineq}~and~\ref{thm:random-walk} in Section~\ref{sec:random-walk}. Some remarks on the structure of almost isoperimetric sets are given in Section~\ref{sec:str}. Theorem~\ref{thm:isop-ineq} is proved in section~\ref{sec:ineq-proof}.

\section{A random walk interpretation of Theorem~\ref{thm:isop-ineq}}
\label{sec:random-walk}

Inequality (\ref{ineq:zero-outside}) is an inequality between two quadratic forms, which can be interpreted as a matrix inequality.

\noi Let $L = L_A$ be the $|A| \times |A|$ matrix indexed by the vertices of $A$, with the following entries: $L(a,a) = n$; and for $a \not = b$, $L(a,b) = -1$ if $a,b$ are connected, and $0$ if not.

\noi Let $J := J_A$ be the $|A| \times |A|$ all-$1$ matrix. Then, (\ref{ineq:zero-outside}) is equivalent to:
\beqn
\label{ineq:Tal-matrix}
L \succeq \frac{1}{|A|} \log_2\(\frac{2^{n}}{|A|}\) \cdot J
\eeqn

\noi This is an inequality of the form $L \succeq v v^t$ for a vector $v \in \R^{A}$. Note that if $A$ is not the complete cube (which we may assume), the matrix $L$ is non-singular. Therefore
$$
L \succeq v v^t ~~~\Leftrightarrow ~~~ I \succeq \(L^{-1/2} v\) \(L^{-1/2} v\)^t  ~~~\Leftrightarrow ~~~
\<L^{-1/2} v, L^{-1/2} v\> \le 1 ~~~\Leftrightarrow ~~~ \<L^{-1} v, v\> \le 1
$$
Let $r = \frac{1}{|A|} \log_2\(\frac{2^{n}}{|A|}\)$ and let $\1$ be the all-$1$ vector in $\R^A$. Then (\ref{ineq:Tal-matrix})  amounts to
\beqn
\label{ident:matrix-form}
\left < L^{-1} \1, \1 \right > \le \frac{1}{r}
\eeqn
This inequality allows a random walk interpretation. Write $L = n\cdot I - E$, where $I$ is the identity matrix and $E$ is the adjacency matrix of a subgraph of $\H$ induced by the vertices in $A$.\footnote{Thus $L$ is the "external" Laplacian of the subgraph induced by $A$.} The matrix $\frac{1}{n} \cdot E$ has eigenvalues smaller than $1$, and therefore we can write
$$
L^{-1} = \frac{1}{n} \cdot \sum_{k=0}^{\infty} \frac{E^k}{n^k}
$$
The inequality (\ref{ident:matrix-form}) can be rewritten as
$$
\frac{n}{r} \ge n \cdot \left < L^{-1} \1, \1 \right > = \sum_{k=0}^{\infty} \frac{\left<E^k \1,\1\right>}{n^k}
$$
Let $Y$ be a random variable defined as follows: choose a uniform random point $a \in A$ and consider a simple random walk in $\H$ starting from $a$. Then $Y$ measures the first time the walk exits $A$. Note that $E^k(a,b)$ counts the number of paths of length $k$ in $A$ between $a$ and $b$. Hence $\frac{1}{n^k} \cdot \sum_{b \in A} E^k(a,b)$ is the probability that the random walk starting from $a$ remains in $A$ for the first $k$ steps, and
$\frac{\left<E^k \1,\1\right>}{|A| \cdot n^k}$ is the probability $Y > k$. Therefore, by (\ref{ident:matrix-form})
$$
\Ex Y = \sum_{k=0}^{\infty} \mbox{Pr}\left\{Y > k\right\} = \frac{1}{|A|} \cdot \sum_{k=0}^{\infty} \frac{\left<E^k \1,\1\right>}{n^k} \le \frac{n}{|A| \cdot r} = \frac{n}{\log_2\(\frac{2^n}{|A|}\)},
$$
proving (\ref{ineq:random-walk}).

\noi Next, we verify that (\ref{ineq:random-walk}) holds with equality if $A$ is a subcube, completing the proof of Theorem~\ref{thm:random-walk}.

\noi Let $A$ be a $d$-dimensional subcube. Then $\mbox{Pr}\left\{Y > k\right\} = \frac{d^k}{n^k}$, and therefore
$$
\Ex Y = \sum_{k=0}^{\infty} \mbox{Pr}\left\{Y > k\right\} = \sum_{k=0}^{\infty} \frac{d^k}{n^k} = \frac{n}{n-d}  = \frac{n}{\log_2\(\frac{2^n}{|A|}\)}
$$
\eprf

\noi One might consider the possibility that subcubes have a stronger property, namely that for a walk of {\it any} length the probability to remain in a subcube is maximal among all sets of the same size. This is true for walks of length $1$, since subcubes have the smallest edge-boundaries. However, the following example shows this to be false already for walks of length $2$:
\xam
\label{xam:bad-2-walk}
The number of length-$2$ walks inside the set $A$ is
\[
\sum_{a,b\in A} E^2(a,b) = \<E^2\1,\1\> = \<E\1,E\1\> = \sum_{x \in A} d^2_x
\]
where $d_x$ is the degree of $x$ in the subgraph induced by $A$. Therefore, for a $d$-dimensional cube, the number of such walks is $2^d \cdot d^2$. But, for a radius-$1$ ball of dimension $2^d - 1$, this number is $\(2^d-1\)^2 + \(2^d - 1\) = 2^d \cdot \(2^d-1\)$, which is much larger.
\exam

\section{Near-isoperimetric sets and their eigenvalues}
\label{sec:str}
Fix a small parameter $\epsilon > 0$. A set $A$ is {\it nearly isoperimetric} if it satisfies the isoperimetric inequality (\ref{ineq:edge-isop}) with nearly an equality, that is
\beqn
\label{ineq:nearly-isop}
|A| \log_2\(\frac{2^{n}}{|A|}\) \le |\partial A| \le (1 + \epsilon) \cdot |A| \log_2\(\frac{2^{n}}{|A|}\)
\eeqn
We would like to understand the structure of nearly-isoperimetric sets and, in particular, their possible similarity to subcubes.

\noi This discussion is closely related to {\it stability} of isoperimetric inequalities. A stability-type result shows that a nearly-isoperimetric set is close (in an appropriate metric) to a genuinely isoperimetric set. Such a result is proved in \cite{Ellis}: Let $\delta$ be at most a small constant, and let $A$ be a set satisfying (\ref{ineq:nearly-isop}) with $\epsilon = \frac{\delta}{\log_2\(2^{n}/|A|\)}$. Then there is a subcube $C$ such that $|A \Delta C| \le O\(\frac{\delta}{ \log\(1/\delta\)} \cdot |A|\)$.

\noi In this section, we look at eigenvalues and eigenvectors of the Laplacian $L$ (equivalently, of the adjacency matrix $E$) of a subgraph induced by an almost isoperimetric subset $A$ of the cube. If $A$ is a subcube, the induced subgraph is regular, of degree $\log_2 |A|$. This means that the minimal eigenvalue of the Laplacian $L$ is $\log_2\(\frac{2^{n}}{|A|}\)$ and the corresponding eigenvector is the all-$1$ vector $\1$.

\noi We show in Corollary~\ref{cor:pseudo-regular} below that if $\epsilon'$ is at most a small constant and $\epsilon = \frac{\epsilon'}{2n} \cdot \log_2\(2^{n}/|A|\)$, then the subgraph induced by a set $A$ satisfying (\ref{ineq:nearly-isop}) is nearly regular, will the degrees of almost all the vertices close to $\log_2 |A|$.

\noi Similar arguments can be used to show that even for $\epsilon$ as large as a small constant, most of the spectral mass in the expansion of $\1$ in an eigenbasis of $L$ is concentrated around the eigenvalue $\log_2\(\frac{2^{n}}{|A|}\)$ (we don't go into details). On the other end of the scale, for a very small $\epsilon \ll \frac{1}{n \cdot \log_2\(\frac{2^{n}}{|A|}\)}$ we can derive stability-type results in the sense of \cite{Ellis} (via a result of Keevash \cite{Keevash} on stability of the Kruskal-Katona inequality). Since this is weaker than the results in \cite{Ellis}, we omit the details here as well.

We start with some notation. Let $|A| = m$, and let $\lambda_1 \ge \lambda_2 \ge \ldots \ge \lambda_m$ be the eigenvalues of $E$. The eigenvalues of $L$ are $n - \lambda_ 1 \le n - \lambda_2 \le \ldots \le n - \lambda_m$. Let $v_1,\ldots, v_m$ be an orthonormal basis of eigenvectors, and let $\1 = \sum_{i=1}^m \alpha_i v_i$ be the expansion of the constant-$1$ function $\1$ in this basis. Note, for future use, that $\sum_{i=1}^m \alpha^2_i = \<\1,\1\> = |A|$.

\noi The inequality (\ref{ident:matrix-form}) translates to
\beqn
\label{ineq:eigen-val}
\sum_{i=1}^m \frac{\alpha^2_i}{n - \lambda_i} \le \frac{1}{r} = \frac{|A|}{\log_2\(\frac{2^n}{|A|}\)}
\eeqn
Note that the edge-boundary of $A$ is given by
\[
|\partial A| = \left<L \1,\1 \right> = \sum_{i=1}^m \alpha^2_i \(n - \lambda_i\)
\]

\noi Therefore, the near-isoperimetric property (\ref{ineq:nearly-isop}) is equivalent to
$$
|A| \log_2\(\frac{2^{n}}{|A|}\) \le \sum_{i=1}^m \alpha^2_i \(n - \lambda_i\)  \le  (1 + \epsilon) \cdot |A| \log_2\(\frac{2^n}{|A|}\)
$$

\noi Consider the probability distribution on $[m]$, given by $p_i = \frac{\alpha^2_i}{|A|}$, and let $f:~i \mapsto n - \lambda_i$ be a positive function on $[m]$. Computing expectations according to $p$, we have $\Ex \frac{1}{f} \cdot \Ex f \le 1 + \epsilon$. Intuitively, this should mean $f$ is concentrated with respect to $p$. In the next lemma we state this formally.

\lem
\label{lem:near-harmonic-inequality}
Let $g$ be a strictly positive-valued function on a finite domain satisfying
$$
\Ex \frac{1}{g} \cdot \Ex g \le 1 + \epsilon
$$
Then
\beqn
\label{ineq-concent-root}
\Ex \(g - \Ex g\)^2 \le \epsilon \cdot \Ex g \cdot \|g\|_{\infty}
\eeqn
\elem
\prf
We have
$$
\Ex\(\frac{(g - \Ex g)^2}{g}\) = \Ex^2 g \cdot \Ex\(\frac{1}{g}\) - \Ex g = \Ex g \cdot \(\Ex g \cdot \Ex\(\frac{1}{g}\) - 1\) \le \epsilon \cdot \Ex g
$$
Therefore
$$
\Ex \(g - \Ex g\)^2 \le  \Ex\(\frac{(g - \Ex g)^2}{g}\) \cdot \|g\|_{\infty} \le \epsilon \cdot \Ex g \cdot \|g\|_{\infty}
$$
\eprf

\cor
\label{cor:pseudo-regular}
Let $A$ satisfy (\ref{ineq:nearly-isop}) with $\epsilon = \frac{\epsilon'}{2n} \cdot \log_2\(\frac{2^n}{|A|}\))$, where $\epsilon' \le 1$. Fix a parameter $0 \le \delta \le 1$.

\noi Choose uniformly at random an element $x \in A$ and consider its outdegree\footnote{the number of neighbors of $x$ outside $A$} $d_{out}(x)$. Then
\[
Pr\left\{(1-\delta) \cdot \log_2\(\frac{2^n}{|A|}\) \le d_{out}(x) \le (1+\epsilon) (1+\delta) \cdot \log_2\(\frac{2^n}{|A|}\)\right\} \ge 1 - \frac{\epsilon'}{\delta^2}
\]
In particular, the subgraph induced by $A$ is almost regular, similarly to the isoperimetric case.
\ecor
\prf
We use the notation above. Consider the random variable $d_{out}(x)$, for $x$ uniformly distributed in $A$. We have
\[
\Ex d_{out}(x) = \frac{|\partial A|}{|A|} = \frac{1}{|A|} \<L\1,\1\>  = \sum_{i=1}^m \frac{\alpha^2_i}{|A|} \(n - \lambda_i\) = \Ex f
\]
Similarly, $\Ex d^2_{out} = \Ex f^2$. Therefore, by Chebyshev's inequality and Lemma~\ref{ineq-concent-root},
\[
Pr\left\{\bigg | d_{out}(x) - \frac{|\partial A|}{|A|} \bigg | \ge \delta \cdot \frac{|\partial A|}{|A|}\right\} = Pr\left\{\Big | d_{out}(x) - \Ex d_{out} \Big | \ge \delta \cdot \Ex d_{out}\right\} \le
\]
\[
\frac{Var\(d_{out}\)}{\delta^2 \cdot \Ex^2 d_{out}} =
\frac{Var\(f\)}{\delta^2 \cdot \Ex^2 f} \le \frac{\epsilon \cdot \|f\|_{\infty}}{\delta^2 \cdot \Ex f} \le \frac{\epsilon'}{\delta^2}
\]
In the last inequality we used the easy fact $\|f\|_{\infty} \le 2n$. The claim follows.
\eprf

\section{Proof of Theorem~\ref{thm:isop-ineq}}
\label{sec:ineq-proof}

\noi There are several simple assumptions we may and will make on the structure of the function $g$ in (\ref{ineq:zero-outside}).

\noi First, we may assume $g \ge 0$, since replacing $g$ with its absolute value preserves RHS of (\ref{ineq:zero-outside}) and can only decrease its LHS.

\noi Second, we may assume the support of $g$ is the whole set $A$, otherwise we may replace $A$ with the support of $g$ in (\ref{ineq:zero-outside}), increasing RHS.

\noi Next, consider the partial order on $\H$ in which $x \preceq y$ iff $x_i \le y_i$, $i = 1,...,n$. A function $g$ on the cube is {\it downwards monotone} if $g(x) \ge g(y)$ when $x \preceq y$.

\noi We may assume the function $g$ in (\ref{ineq:zero-outside}) to be monotone. This follows from two simple lemmas.

\lem
\label{lem:one-shift}
Fix a direction $1 \le i \le n$, and let $f$ be a function obtained from $g$ by a downward shift in direction $i$.

\noi That is, for any pair of adjacent points $x,y$ in the cube, with $x_i = 0$ and $y_i = 1$, set
\[
f(x) = \max\{g(x),g(y)\}~~~\verb+and+~~~f(y) = \min\{g(x),g(y)\}
\]

\noi Then
\[
\E(f,f) \le \E(g,g)
\]
\elem
\prf
This is a standard "shifting" argument \cite{Bollobas}, more commonly applied in the special case of $g$ being a characteristic function. The claim of the lemma is easily seen to follow from its validity for $2$-dimensional cubes. The two-dimensional case is verifiable by a direct calculation.
\eprf

\lem
\label{lem:many-shifts}
Applying consecutive shifts in directions $i = 1,...,n$ to a function on the cube produces a monotone function.
\elem
\prf
Again, it suffices to verify this in the two-dimensional case. See \cite{mon-test} where this argument is applied in the special case of characteristic functions.
\eprf

\noi The proof proceeds by induction on the dimension.

\noi First, consider the base case $n=1$. There are two choices for $|A|$. If $|A| = 1$, we are in the boolean case, in which (\ref{ineq:zero-outside}) is the usual edge-isoperimetry. If $|A| = 2$, RHS in (\ref{ineq:zero-outside}) is $0$, and we are done.

\noi Now we go to the induction step.

\noi The cube $\H$ decomposes into two $(n-1)$-dimensional subcubes. The first subcube contains all vectors with last coordinate $0$, and the second all vectors with last coordinate $1$. The function $g$ and the set $A$ decompose according to their restrictions to the subcubes.
$$
g \hookrightarrow \(g_0, g_1\),~~~A \hookrightarrow \(A_0, A_1\)
$$

\noi Induction step amounts to proving
$$
\E(g,g) = \frac12 \cdot \(\E\(g_0,g_0\) + \E\(g_1,g_1\)\) + \|g_0 - g_1\|^2_2 \ge_{{\rm ind}}
$$
$$
\frac12 \cdot 2 \cdot \frac{1}{2^{n-1} |A_0|} \log\(\frac{2^{n-1}}{|A_0|}\) \(\sum_{x \in A_0} g_0(x)\)^2 +
$$
$$
\frac12 \cdot
2 \cdot \frac{1}{2^{n-1} |A_1|} \log\(\frac{2^{n-1}}{|A_1|}\) \(\sum_{x \in A_1} g_1(x)\)^2 + \|g_0 - g_1\|^2_2 \ge^{{\bf ??}}
$$
$$
2 \cdot \frac{1}{2^n |A|} \log\(\frac{2^n}{|A|}\) \(\sum_{x \in A} g(x)\)^2
$$
In the expressions above, the Dirichlet forms and the $\ell_2$ distance for functions $g_i$ on $(n-1)$-dimensional cubes are computed with respect to the uniform probability measure on these subcubes.

\noi Note that, by our assumptions on $g$, the set $A$ is downwards monotone, since it is the support of a monotone function $g$. This implies $A_1 \subseteq A_0$ (identifying the two subcubes in the natural way).

\noi The expression we need to analyze allows an additional simplifying assumption on $g$. We may assume both $g_0, g_1$ to be constant on $A_1$ and on $A_0 \setminus A_1$ (and of course $g_i$ vanish on $A^c_i$, in particular $g_1$ is zero on $A_0 \setminus A_1$). In fact, replacing $g_i$ with their averages on the corresponding subsets can only decrease LHS and does not change RHS in the second inequality above.

\noi We proceed with analysis, introducing some notation.

\noindent {\bf Notation}:
\begin{itemize}
\item
Let $s_0 := \sum_{x \in A_0} g_0(x)$, $s_1 := \sum_{x \in A_1} g_1(x)$. Let $t_0 := |A_0|$, $t_1 := |A_1|$. We may and will assume $t_1 > 0$ and $s_1 > 0$, otherwise the problem reduces to a lower-dimensional case.
\item
Let $\alpha$ be the value of $g_0$ on $A_1$ and $\gamma$ be the value of $g_0$ on $A_0 \setminus A_1$. Let $\beta$ be the value of $g_1$ on $A_1$.
\item
The "$f$"-notation. Let $f(t) = f_{n-1}(t):= \frac{1}{t} \log_2\(\frac{2^{n-1}}{t}\)$.
\end{itemize}

\noi Note that
\begin{enumerate}
\item
$$
t_1 \beta = s_1
$$
\item
$$
t_1 \alpha + \(t_0 - t_1\) \gamma = s_0
$$
\item
$$
\|g_0 - g_1\|_2^2 = \frac{1}{2^{n-1}} \cdot \(t_1 \(\alpha-\beta\)^2 + \(t_0 - t_1\) \gamma^2\)
$$
\end{enumerate}

\noi With the new notation, the inequality to be verified for the induction step is:

\beqn
\label{ineq:ind-step-f}
f\(t_0\) s^2_0 + f\(t_1\) s^2_1 + \(t_1 \(\alpha-\beta\)^2 + \(t_0 - t_1\) \gamma^2\) \ge \frac12 \cdot f\(\frac{t_0 + t_1}{2}\) \(s_0 + s_1\)^2
\eeqn

\noi Expressing $\beta$ and $\gamma$ as functions of $s_i, t_i$ and of $\alpha$,\footnote{Dealing with the simple case $t_0 = t_1$ separately} LHS of (\ref{ineq:ind-step-f}) is a quadratic in $\alpha$ with coefficients depending on $s_i$ and $t_i$. Minimizing LHS in $\alpha$, we arrive, after some simple calculations, to the following inequality we need to verify:
\beqn
\label{ineq:t-constraints}
f\(t_0\) s^2_0 + f\(t_1\) s^2_1 + \frac{\(s_0-s_1\)^2}{t_0} \ge \frac12 \cdot f\(\frac{t_0 + t_1}{2}\) \(s_0 + s_1\)^2
\eeqn

\noi Next, let $R = \frac{s_0}{s_1}$.

\noi Inequality (\ref{ineq:t-constraints}) transforms to a quadratic inequality in $R$
\beqn
\label{R-quad}
f\(t_0\) R^2 + f\(t_1\) + \frac{(R-1)^2}{t_0} \ge \frac12 \cdot f\(\frac{t_0 + t_1}{2}\) \(R + 1\)^2
\eeqn
We need to check $P(R) := aR^2 + bR + c \ge 0$ with the coefficients $a,b,c$ coming from (\ref{R-quad}). We will, in fact, verify $a \ge 0$ and $D = b^2 - 4ac \le 0$, which will conclude the proof.

\noi We start with some simple properties of the function $f(t) = \frac{1}{t} \log_2\(\frac{2^{n-1}}{t}\)$.
\lem
\label{lem:f-prop}
The function $f(t)$ is decreasing and convex for $0 < t < 2^{n-1}$. It satisfies the identity
\beqn
\label{iden:f}
f(\beta \cdot t) = \frac{1}{\beta} \cdot f(t) + \frac{1}{\beta} \log \frac{1}{\beta} \cdot \frac{1}{t}
\eeqn
for any $t, \beta > 0$.
\elem
\prf
Directly verifiable.
\eprf

\cor
\label{cor:infty}
Viewing inequality (\ref{R-quad}) in the form $aR^2 + bR + c \ge 0$, we have
$$
a \ge 0
$$
\ecor
\prf
It is easy to verify
$$
a = \frac{2 t_0 f\(t_0\) + 2 - t_0 f\(\frac{t_0 + t_1}{2}\)}{2t_0}
$$
By Lemma~\ref{lem:f-prop}
\[
t_0 f\(\frac{t_0 + t_1}{2}\) \le t_0 f\(\frac{t_0}{2}\) = 2 t_0 f\(t_0\) + 2
\]
completing the proof.
\eprf

\noi It remains to verify the inequality $4ac \ge b^2$, which, after some simplification, reduces to:
\beqn
\label{ineq:technical-main}
2t_0 f\(t_0\) f\(t_1\) + 2(f\(t_0\) + f\(t_1\)) \ge t_0 f\(\frac{t_0 + t_1}{2}\)(f\(t_0\) + f\(t_1\)) + 4 f\(\frac{t_0 + t_1}{2}\)
\eeqn

\section{Proof of inequality (\ref{ineq:technical-main})}
Renaming the variables $x = t_0$ and $y = t_1$, and recalling the constraints on $t_0$ and $t_1$, we need to prove (\ref{ineq:technical-main}) for $1 \le y < x \le 2^{n-1}$.

\noi Rearranging, this is easily seen to be equivalent to
\beqn
\label{ineq: with Delta}
\Delta(x,y) \ge \frac{x\cdot\(f(x) - f(y)\)^2}{2x\cdot\(f(x) + f(y)\) + 8}
\eeqn
Here $\Delta(x,y) := \frac{f(x) + f(y)}{2} - f\(\frac{x+y}{2}\)$. Note that $\Delta \ge 0$ since $f$ is convex.

\noi We now substitute $y = \beta x$ in (\ref{ineq: with Delta}), witn $0 < \beta < 1$ and expand using (\ref{iden:f}). We have
$$
\Delta(x,y) = \Delta(x,\beta \cdot x) = \frac{f(x) + f(\beta\cdot x)}{2} - f\(\frac{1+\beta}{2} \cdot x\) =
$$
$$
\frac12 \cdot \(f(x) + \frac{1}{\beta} \cdot f(x) + \frac{1}{\beta} \log \frac{1}{\beta} \cdot \frac{1}{x}\) - \(\frac{2}{1+\beta} \cdot f(x) + \frac{2}{1+\beta} \log \frac{2}{1+\beta} \cdot \frac{1}{x}\) =
$$
$$
\frac{(1-\beta)^2}{2\beta(1+\beta)} \cdot f(x) + \(\frac{1}{2\beta}\log \frac{1}{\beta} - \frac{2}{1+\beta} \log \frac{2}{1+\beta} \)\cdot \frac{1}{x}
$$
As to RHS of (\ref{ineq: with Delta}), we have
$$
RHS(x,y) = RHS(x,\beta \cdot x) = \frac{x\cdot \(\frac{1}{\beta} \cdot f(x) + \frac{1}{\beta} \log \frac{1}{\beta} \cdot \frac{1}{x} -f(x)\)^2}{2x \cdot\(f(x) + \frac{1}{\beta} \cdot f(x) + \frac{1}{\beta} \log \frac{1}{\beta} \cdot \frac{1}{x}\) + 8}
$$
Taking $z := xf(x)$,
$$
\Delta \ge RHS ~~~ \Longleftrightarrow ~~~ x \cdot \Delta \ge x \cdot RHS ~~~ \Longleftrightarrow ~~~ A z + B \ge \frac{(Cz +D)^2}{Ez+F},
$$
where $z = xf(x) = \log\(\frac{2^{n-1}}{x}\) \ge 0$ and $A,B,...,F$ depend only on $\beta$.

Specifically,
$$
\left\{ \begin{array}{lll}
A & = & \frac{(1-\beta)^2}{2\beta(1+\beta)} \\
B & = & \frac{1}{2\beta}\log\frac{1}{\beta} - \frac{2}{1+\beta} \log \frac{2}{1+\beta}\\
C & = & \frac{1-\beta}{\beta}\\
D & = & \frac{1}{\beta} \log \frac{1}{\beta}\\
E & = & \frac{2 + 2\beta}{\beta}\\
F & = & \frac{2}{\beta} \log \frac{1}{\beta} + 8\\
\end{array} \right.
$$
So, we need to show
\[
(Az + B)(Ez+F) \ge (Cz+D)^2
\]
Observe that $AE = C^2 = \frac{(1-\beta)^2}{\beta^2}$. Therefore, this reduces to a linear inequality in $z$:
\[
(AF + BE - 2CD)\cdot z\ge D^2 - BF
\]
This holds for all nonnegative $z$ if and only if
$$
\left\{ \begin{array}{lll}
AF + BE & \ge & 2CD\\
BF & \ge & D^2\\
\end{array}\right.
$$
Hence, the problem is reduced to two univariate inequalities in $\beta$. We will prove them in the next two lemmas.
\lem
$$
AF + BE \ge 2CD
$$
For $0 < \beta < 1$.
\elem
\prf
Simplifying and rearranging, this inequality reduces to
$$
\frac{(1-\beta)^2}{\beta(1+\beta)} \cdot \log \frac{1}{\beta} + \frac{1+\beta}{\beta} \cdot \log \frac{1}{\beta} + 4 \frac{(1-\beta)^2}{1 + \beta} \ge 4\log\frac{2}{1+\beta} + 2\frac{1-\beta}{\beta} \cdot \log \frac{1}{\beta}
$$
$$
\frac{\beta}{1+\beta} \cdot \log\frac{1}{\beta} + \frac{(1-\beta)^2}{1 + \beta} \ge \log\frac{2}{1+\beta}
$$
$$
\beta \log \frac{1}{\beta} + (1-\beta)^2 \ge (1+\beta) \log\frac{2}{1+\beta}
$$
The derivative of $g(\beta) = \beta \log \frac{1}{\beta} + (1-\beta)^2 - (1+\beta) \log\frac{2}{1+\beta}$ is $\log\frac{1+\beta}{2\beta} - 2(1-\beta)$. This is a convex function, which means it can vanish in at most two points in the interval $(0,1]$. In addition, $g'$ is positive close to $0$ and it vanishes at $1$. Taking into account the boundary conditions $g(0) = g(1) = 0$, this means that $g$ first increases from $0$ at zero and then decreases to $0$ at one, that is, it is nonnegative.

\eprf

\lem
$$
BF \ge D^2
$$
for $0 < \beta < 1$.
\elem
\prf
We need to prove
$$
\(\frac{1}{2\beta}\log\frac{1}{\beta} - \frac{2}{1+\beta} \log \frac{2}{1+\beta}\)\cdot\(\frac{2}{\beta} \log \frac{1}{\beta} + 8\) \ge \frac{1}{\beta^2} \log^2\frac{1}{\beta}
$$
Simplifying and rearranging, this reduces to
$$
(1+\beta) \cdot \log\frac{1}{\beta} \ge \log\frac{2}{1+\beta} \cdot \log\frac{1}{\beta} + 4\beta\cdot \log\frac{2}{1+\beta}
$$
$$
\(\beta + \log(1+\beta)\)\cdot \log\frac{1}{\beta} \ge 4 \beta \cdot \log\frac{2}{1+\beta}
$$
As in the preceding lemma, the function $g(\beta) = \(\beta + \log(1+\beta)\)\cdot \log\frac{1}{\beta} - 4 \beta \cdot \log\frac{2}{1+\beta}$ vanishes at the endpoints. We will (again) claim it increases from $0$ at zero and then decreases from the maximum point to $0$ at one, and is, therefore, nonnegative on the interval.

\noi As before, it will suffice to show that $g'$ is convex, is positive in the beginning of the interval, and vanishes at $1$. We have
$$
\ln 2 \cdot g'(\beta) = \(1+\frac{1}{\ln 2\cdot(1+\beta)}\) \cdot \ln\frac{1}{\beta} - \frac{\ln 2 \cdot \beta + \ln(1+\beta)}{\ln 2 \cdot \beta} - 4\ln \frac{2}{1+\beta} + \frac{4\beta}{1+\beta}
$$
It is easy to verify that $g'$ is positive for small positive $\beta$ and that $g'(1) = 0$.

\noi It remains to check $g'$ is convex. Taking another two derivatives, we have
$$
\ln 2 \cdot g'''(\beta) = \(\frac{1}{\beta^2} + \frac{3}{\ln 2} \cdot \(\frac{1}{\beta^2(1+\beta)} + \frac{1}{\beta(1+\beta)^2}\) + \frac{2}{\ln 2} \cdot \frac{\ln \frac{1}{\beta}}{(1+\beta)^3}\) -
$$
$$
\(\frac{4}{(1+\beta)^2} + \frac{8}{(1+\beta)^3} + \frac{2}{\ln 2} \cdot \frac{\ln (1+\beta)}{\beta^3}\)
$$
To show this is nonnegative, we multiply by $\beta^3 (1+\beta)^3$ and verify
$$
\beta(1+\beta)^3 + \frac{3}{\ln 2} \cdot \beta (1+\beta)(1+2\beta) + \frac{2}{\ln 2} \cdot \beta^3 \ln \frac{1}{\beta} \ge
4 \beta^3 (1+\beta) + 8 \beta^3 + \frac{2}{\ln 2} \cdot (1+\beta)^3 \ln (1+\beta)
$$
We show a stronger inequality (removing third summand on the left)
$$
\beta(1+\beta)^3 + \frac{3}{\ln 2} \cdot \beta (1+\beta)(1+2\beta)  \ge
4 \beta^3 (1+\beta) + 8 \beta^3 + \frac{2}{\ln 2} \cdot (1+\beta)^3 \ln (1+\beta)
$$
Note that $(1+\beta)\ln(1+\beta) \le 2\ln 2 \cdot \beta$, by convexity of $(1+\beta)\ln(1+\beta)$ on $[0,1]$. Substituting and simplifying, it suffices to show
$$
(1+\beta)^3 + \frac{3}{\ln 2} \cdot (1+\beta)(1+2\beta)  \ge
4 \beta^2 (1+\beta) + 8 \beta^2 + 4 \cdot (1+\beta)^2
$$
Since $\frac{3}{\ln 2} \ge 4$ and $\beta^2 \ge \beta^3$ for $0 \le \beta \le 1$, it suffices to prove the quadratic inequality
$$
1 + 3\beta + 3\beta^2 + 4 (1+\beta)(1+2\beta) \ge 15 \beta^2 + 4 \cdot (1+\beta)^2
$$
Simplifying, this reduces to the trivial statement
$$
7\beta + 1 \ge 8 \beta^2
$$
\eprf

\end{document}